\documentclass[12pt,a4paper,oneside,reqno]{article}

\usepackage{amsmath,amssymb}
\usepackage{amsthm}
\usepackage{color}
\usepackage{endnotes}
\usepackage[xdvi]{graphicx}
\usepackage{psfrag}

\newcommand{\ax}[1]{\textcolor{axcolor}{\ensuremath{\mathsf{#1}}}}

\definecolor{axcolor}{rgb}{.4,0,.4}
\definecolor{axbgcolor}{rgb}{1,.7,1}

\definecolor{bl}{rgb}{0,0,1}
\definecolor{gr}{rgb}{0,.5,0}
\newcommand{\AMbl}{\textcolor{bl}}

\newcommand{\R}{\mathbb{R}}

\newcommand{\Ph}{\mathsf{Ph}}
\newcommand{\Q}{{Q}}
\newcommand{\B}{{B}}
\newcommand{\W}{\mathsf{W}}
\newcommand{\IB}{\mathsf{IB}}
\newcommand{\IOb}{\mathsf{IOb}}
\newcommand{\Ob}{\mathsf{Ob}}

\newcommand{\Setclose}{\,\right\}}

\newcommand{\Setopen}{\left\{\,}
\newcommand{\leteq}{\mbox{$:=$}}

\newcommand{\then}{\rightarrow}
\renewcommand{\iff}{\leftrightarrow }
\newcommand{\defiff}{\stackrel{def}{\Longleftrightarrow}}

\newcommand{\comment}[1]{}

\theoremstyle{definition} \newtheorem{thm}{Theorem}%[section]
\theoremstyle{definition} %[section]
\theoremstyle{definition} %[section]

\begin{document}

\makeatletter
 \renewcommand\section{\@startsection {section}{1}{\z@}%
                     {-3.25ex\@plus -1ex \@minus -.2ex}%
                     {1.5ex \@plus .2ex}%
                     {\normalfont\scshape\centering}}
\makeatother

\title{Vienna Circle and \\ Logical Analysis of Relativity Theory}
%Deriving Axioms of General Relativity from that of Special Relativity}
\author{H.~Andr{\'e}ka, J.~X.~Madar{\'a}sz, I.~N{\'e}meti, P.~N{\'e}meti and G.~Sz{\'e}kely}
\date{}

\maketitle

\section{introduction}

In this paper we present some of our school's results in the area of
building up relativity theory (RT) as a hierarchy of theories in the
sense of logic. We use plain first-order logic (FOL) as in the
foundation of mathematics (FOM) and we build on experience gained in
FOM.
The main aims of our school are the following:  We want to base the
theory on simple, unambiguous axioms with clear meanings. It should
be absolutely understandable for any reader what the axioms say and
the reader can decide about each axiom whether he likes it. The
theory should be built up from these axioms in a straightforward,
logical manner. We want to provide an analysis of the logical
structure of the theory. We investigate which axioms are needed for
which predictions of RT. We want to make RT more transparent
logically, easier to understand, easier to change, modular, and
easier to teach. We want to obtain deeper understanding of RT.

Our work can be considered as a case-study showing that the Vienna
Circle's (VC) approach to doing science is workable and fruitful
when performed with using the insights and tools of mathematical
logic acquired since its formation years at the very time of the VC
activity. We think that logical positivism was based on the insight
and anticipation of what mathematical logic is capable when
elaborated to some depth. Logical positivism, in great part
represented by VC, influenced and  took part in the birth of modern
mathematical logic. The members of VC were brave forerunners and
pioneers.

Let's see what was available before or during the VC activity, and
what was not available for the members of VC but available for us
now.
The VC activities in the strict sense fall into the period of
1922-1936. The beginning of intensive development of FOL coincides
with this period.
The elements of first-order language FOL (propositional logic,
quantifiers) were available based on the works of Boole, Peirce,
Schr\"oder, Frege, Russel (roughly 1860-1910). However, the
completeness and incompleteness theorems, compactness theorem,
semantics for FOL, model theory, proof theory and definability
theory were not there before VC. Many of these became available as
works of people influenced by or working in VC. Here is a brief
chronological order for some turning-points in this development:
1929 G\"odel's completeness theorem, 1930 Tarski's decision method
for the elementary theory of reals, 1931 G\"odel's incompleteness
theorem, 1933 Tarski mathematical definition of truth, 1936 Tarski
concept of semantical consequence relation, definition of model
theoretic semantics. Proof theory developed only later than around
1940, and model theory developed only later than around 1950. Beth's
definability theorem was anticipated by %VC associate
Hans Reichenbach 1924 (motivated by RT), Tarski 1936, but became
available only in 1952. The authors of the standard Model Theory
book, Chang and Keisler, are both students of Tarski. Two chapters
of model theory, the theory of semantics and the theory of
definability are of great importance for our work.
So what was/is FOL used for? Some dates showing the emergence of
paradoxes in mathematics necessitating an axiomatic approach,
experiments with formal languages to describe, found and unify
mathematics and then all science are as follows: 1897 Burali-Forti
Paradox, 1900 Hilbert's Program, 1901 Russel's Paradox, 1908
Zermelo's Axiom system for Set Theory, 1910 Russel-Whitehead
Principia Mathematica, 1918 L\"owenhein-Skolem ``Paradox," 1922
Fraenkel's addition to Zermelo-Fraenkel Set theory, 1924
Tarski-Banach ``Paradox,"  1926 Tarski's axiomatization of Euclidean
geometry in FOL, 1930 the addition of the Axiom of Regularity to Set
Theory, beginning with 1935 the Bourbaki group's work formalizing
and uniting mathematics, 1937 Tarski role of logic in scientific
studies. Success story of FOM. G\"odel-Bernays and von
Neumann's Set Theory. %(evszamok, utana)
We will see in this paper many of the above results used in our
work. G\"odel, Tarski, Reichenbach, Hilbert, Russel, and Einstein
were all connected to VC in some way or other.

Why relativity theory? Logical positivism is a philosophy which
holds that the only authentic knowledge is that based on
observation, experience, experiment through formal logic. Einstein's
relativity theory transforms space and time from being a priori
Euclidean and absolute things to something which emerges from
experience and experiments, thus claiming the subjects of space and
time from the realm of metaphysics for science. Einstein's
thought-experiments served to bring logic into the picture.
Formalizing relativity theory in logical language was a primary
interest in logical positivism, see Reichenbach's works. Relativity
theory also led to Modern High-Precision Cosmology as a branch of
hard-core physics (and not part of metaphysics). Since space-time is
the arena in which most processes studied by modern science unfold,
a logic based foundation for RT (the theory of space-time) might be
a natural starting point for a foundation and unification of the
whole of science (a VC goal \cite{uni1}, \cite{uni2}).

Our group investigates a hierarchy of relativity theories, weaker
and stronger theories. We not only give axiom systems and prove
their completeness with respect to their intended models, but we
also derive RT's main predictions, ask ourselves which axioms play
the key role in their derivations, we make ``reverse relativity" in
analogy with ``reverse mathematics," and we analyze the theories in
many ways. In this paper we present three of our main axiom systems
(i.e., theories) with just stating some of their most important
properties. These three theories have the same language. We begin
with introducing this language.

\comment{There are great many axiomatizations of Special Relativity
(SR) in the literature, see, e.g., \cite{Guts}, \cite{Mundy},
\cite{Robb}, \cite{Schutz73}, \cite{Schutz}, \cite{Suppes} some of
them are formulated within first-order logic (FOL), see \cite{Ax},
\cite{Goldblatt}, \cite{Pambuccian}... However, there are only a few
axiomatizations of General Relativity (GR), see, e.g.,
\cite{Mould}..., and a very few of them (besides ours) are
formulated within FOL, see, e.g., \cite{Basri}, \cite{Benda}.
Moreover, as far as we know, no paper in the literature deals with
relations of the axioms of the two theories. Here we are going to
fill this gap by deriving the axioms of GR of \cite{logst},
\cite{myphd} from the axioms of SR of \cite{AMNsamples} in two
simple steps.  Both of these steps are natural and well motivated by
some central ideas of Einstein's. In the first step we are going to
extend our axiomatization of SR to accelerated observers. This step
provides us a theory which nicely fills the gap between SR and GR.
In the second step we are going to arrive to a first-order theory of
GR by eliminating the distinction between inertial and noninertial
observers in the level of axioms.
Another aim of this paper is showing that our axiomatization of SR
and GR is not just complete with respect to the intended models, but
also suitable for analyzing the logical structure of these theories,
i.e., analyzing logical connections between the axioms and their
predictions. }%endcomment

\section{the common language of the theories presented here}

We will use FOL. FOL can be viewed as a fragment of natural language
with unambiguous syntax and semantics. One of acknowledged traits of
using FOL is that it helps eliminate tacit assumptions, one of VC's
maxim. The most important decision in writing up an axiom system in
FOL is to choose the vocabulary, or primitive symbols of our
language, i.e., what objects and what relations between them will
belong to the language we will use.

We want to talk about space and time as relativity theory conceives
them. We will talk about space-time as experienced through motion.
We represent motion as changing spatial location in time. We will
call the entities that do the motion ``test-particles." Sometimes,
for using a shorter word, we will call them ``bodies" but in reality
they can be anything that move, e.g., they can be coordinate systems
or electromagnetic waves, or light signals or centers of mass.%
\footnote{Note about extended bodies: We concentrate
on test-particles and regard test-particles as spatially point-like,
i.e., of size zero. As far as we are aware of it, this idealization
is harmless from the point of view of the goals of relativity
theory. If we want to treat an extended body in our theory (as we do
in the theory \ax{AccRel} of accelerated observers), we represent it
as a ``cloud" of test-particles. This is consistent with the spirit
of standard physical worldview of regarding extended bodies as
clouds of elementary particles.}
 To talk about spatial locations and time we will use quantities arranged
 in a (space-time) coordinate system, and we will have a basic relation, the so-called
 worldview relation, which tells us which test-particles are present in which
 locations at which instants. We will think of the quantities as
 the real numbers (i.e., the number-line), so we will use a ``less than" relation and two operations,
 addition and  multiplication, on them. In this paper to axiomatize special
 relativity theory, we will use two more primitive notions, namely that of ``inertial
 test-particles" and ``light-signals" which we will simply call photons.%
\footnote{To talk about light-signals is not necessary for building up
SR. One simple way of avoiding them is defining light-signal as
anything that moves with ``speed of light." There are deeper ways of
avoiding the use of light-signals in building up relativity theory,
see, e.g., \cite[sec.5]{pezsgo}.}

To concretize what we said so far, let us consider the following
two-sorted first-order language:
\[\{\, \B, \IB, \Ph, \Q,+,\cdot,<, \W\,\},\]
 where $\B$ (test-particles or bodies) and $\Q$ (quantities) are the two sorts,
 $\IB$ (inertial bodies) and
 $\Ph$ (light signals or photons) are unary relation symbols of sort
 $\B$,  $\cdot$ and $+$ are binary function symbols and $<$ is a binary
 relation symbol of sort $\Q$,  and $\W$ (the worldview relation) is a $6$-ary
 relation symbol of sort $\B\B\Q\Q\Q\Q$. $\B$ and $\Q$ can be thought
 of as the physical and as the mathematical universes.

 Atomic formulas $\IB(c)$ and $\Ph(p)$  are translated as
 ``$c$ is an inertial body,'' and  ``$p$ is a photon,''
 respectively.  We use the worldview
 relation $\W$ to speak about coordinatization by translating
 $\W(o,b,x,y,z,t)$ as ``observer $o$ coordinatizes body $b$ at
 space-time location $\langle x,y,z,t\rangle$,'' (i.e., at space
 location $\langle x,y,z\rangle$ and at instant $t$). We sometimes use the more intuitive
 expressions ``sees" or ``observes" for coordinatizes. We will use
 the letters, and their variants, $o,b,p,m,k$ for variables of sort
 $\B$, and the letters $x,y,z,t$ and their variants for variables of
 sort $\Q$. For easier readability, we will use $\bar x, \bar y$ for
 sequences of four variables $x_1,x_2,x_3,x_4$ and
 $y_1,y_2,y_3,y_4$.

We have not introduced the concept of observers as a basic one
because it can be defined as follows: an {\bf observer} is nothing
else than a body who ``observes" (coordinatizes) some other bodies
somewhere, this property can be captured by the following
first-order formula of our language:
\begin{equation*}
\Ob(o) \defiff  \exists b\bar x\enskip \W(o,b,\bar x);
%\Ob(o) \defiff \exists bxyzt \enskip W(o,b,x,y,z,t);
\end{equation*}
and {\bf inertial observers} can be defined as inertial bodies which are
observers, formally:
\begin{equation*}
\IOb(o) \defiff \IB(o)\land \Ob(o).
\end{equation*}

 To abbreviate formulas of FOL we often
omit parentheses according to the following convention. Quantifiers
bind as long as they can, and $\land$ binds stronger than
$\rightarrow$. For example, we write $\forall x\enskip
\varphi\land\psi\rightarrow\exists y\enskip \delta\land\eta$ instead
of $\forall x\big((\varphi\land\psi)\rightarrow\exists
y(\delta\land\eta)\big)$.

\section{axioms of special relativity}

Having specified the language, let us turn to the axioms of our
first theory. This will be an axiom system for Special Relativity
theory (SR).

%The first axiom, \ax{AxFrame}, comprises the most basic intuitions
%mentioned so far. We put them into an explicit axiom because we want
%to avoid the use of \mbox{\it any} tacit assumptions, which is a
%Vienna Circle maxim that we follow in all our works.

\begin{description}
\item[\ax{AxField}:]  The
  {\it quantity part} $\langle \Q;+,\cdot, < \rangle$ is an ordered
  field.
\end{description}

\noindent For the FOL definition of linearly ordered field see,
e.g., \cite[p.41]{chang-keisler}; this is a formulation of some of
the most basic properties of addition and multiplication of real
numbers. One of these properties is that there is a unique neutral
element for addition ($\exists z\forall x\, z+x=x$), we call this
element $z$  zero and we denote it with $0$.

The next axiom simply states that each inertial observer assumes
that it rests at the origin of the space part of its coordinate
system. It also can be thought of as expressing that we identify a
coordinate system (or reference frame) with a test-particle
``sitting" at the origin.

\begin{description}
\item[\ax{AxSelf}:] Any inertial observer coordinatizes (observes) itself as
``living on the time-axis," i.e., it coordinatizes itself at a
  coordinate point if and only if the space component of this point is the
  origin:%\\
  %,  i.e., space location $\langle 0,0,0\rangle$:
\begin{equation*}
\forall oxyzt\enskip  \IOb(o)\then \big(\W(o,o,x,y,z,t)\iff
x=y=z=0\big).
\end{equation*}
\end{description}

Our next axiom is on the constancy of the speed of light. For
convenience, we choose $1$ for this speed. This choice physically
means using units of distance compatible with units of time, such as
light-year, light-second, etc.
\begin{description}
\item[\ax{AxPh}:] The speed of light signals is $1$
 and it is possible to ``send out" a photon in any direction, according to any
  inertial observer:
\begin{multline*}
\forall o\bar x\bar x'\enskip \IOb(o)\then \big(
 \exists p(\Ph(p)\land \W(o,p,\bar x)
    \land \W(o,p,\bar x'))\\\iff
  (x_1-x'_1)^2+(x_2-x'_2)^2+(x_3-x'_3)^2=(x_4-x'_4)^2\big).
\end{multline*}
\end{description}
 This is the most important axiom of SR, it is its ``physical" axiom.
 Axiom \ax{AxPh} is very well confirmed by experiments, such as the
 Michaelson--Morley experiment and its variants.
The next axiom establishes connections between different coordinate
systems. It expresses the idea that all observers ``observe" the
same outside reality.
\begin{description}
\item[\ax{AxEv}:] All inertial observers coordinatize the
  same ``meetings of bodies:"
\begin{equation*}
\forall oo'\bar x\enskip \IOb(o)\land\IOb(o')\then \exists
  \bar x'\enskip \forall b \enskip \W(o,b,\bar x)\iff
  \W(o',b,\bar x').
%  $\forall oo'xyzt\enskip \exists
%  x'y'z't'\enskip \forall b \enskip W(o,b,x,y,z,t)\iff
%  W(o',b,x',y',z',t').$
\end{equation*}
\end{description}

We call ``meetings of bodies" events. By our next axiom, we assume
that inertial observers use the same units of measurement. This is
only a ``simplifying" axiom.
\begin{description}
\item[\ax{AxSymd}:] {Inertial} observers agree as to the spatial
  distance between events if these events are simultaneous for both of
  them, formally:
\begin{multline*}
\forall oo'\bar x\bar x'\bar y\bar y'\enskip \IOb(o)\land\IOb(o')\land x_4=y_4\land x'_4=y'_4\land\\
 \forall b \enskip (\W(o,b,\bar x)\iff
    \W(o',b,\bar x'))\land
  \forall b \enskip (\W(o,b,\bar y)\iff
    \W(o',b,\bar y'))\\ \then
(x_1-y_1)^2+(x_2-y_2)^2+(x_3-y_3)^2=(x'_1-x'_2)^2+(y'_1-y'_2)^2+(z'_1-z'_2)^2.
%\forall
%oo'tt'x_1y_1z_1x_2y_2z_2x'_1y'_1z'_1x'_2y'_2z'_2\enskip\\
% \forall b \enskip (W(o,b,x_1,y_1,z_1,t)\iff
%    W(o',b,x'_1,y'_1,z'_1,t'))\land\\
%  \forall b \enskip (W(o,b,x_2,y_2,z_2,t)\iff
%    W(o',b,x'_2,y'_2,z'_2,t'))\\ \then
%(x_1-x_2)^2+(y_1-y_2)^2+(z_1-z_2)^2=(x'_1-x'_2)^2+(y'_1-y'_2)^2+(z'_1-z'_2)^2.
\end{multline*}

\end{description}

Let us now introduce our axiom system of SR as the set of the axioms
above:
\begin{equation*}
\boxed{\ax{SpecRel}=\Setopen \ax{AxField}, \ax{AxSelf}, \ax{AxPh},
  \ax{AxEv},\ax{AxSymd}\Setclose}.
\end{equation*}

The reader is invited to check that all the axioms of \ax{SpecRel}
are simple, comprehensible and observationally oriented. In setting
up an axiom system, we want the axioms be streamlined, economical,
transparent and few in number. On the other hand, we want to have
all the surprising, shocking, paradoxical predictions of RT as
theorems (and not as axioms). We want to have the price-value ratio
to be good, where the axioms are on the ``cost"-side, and the
theorems are on the ``gain"-side.

Let us see what theorems we can prove from \ax{SpecRel}. We will see
that we can prove everything from our five axioms that ``usual" SR
can, but let us proceed more slowly. In the axioms we did not
require explicitly, but it can be proved from \ax{SpecRel} with the
rigorous methods of FOL that inertial observers see each other move
on a straight line, uniformly (covering the same amount of distance
in the same amount of time). For a ``fancy theorem" from ``plain
axioms," let us prove from \ax{SpecRel} that ``no inertial observer
can move faster than light." Below, $\vdash$ denotes derivability in
one of FOL's standard proof systems.

\begin{thm}(\mbox{\sf NoFTL}) In an inertial observer $m$'s
worldview, any inertial observer $k$ moves slower than any
light-signal $p$, i.e., if both $k$ and $p$ move from spatial
locations $\langle x_1,x_2,x_3\rangle$ to $\langle
y_1,y_2,y_3\rangle$, then for the observer $k$ this trip took more
time than for the photon $p$. Formally:
\begin{multline*}
\ax{SpecRel}\vdash
\forall mkp\bar x\bar yt\enskip \IOb(m)\land\IOb(k)\land\Ph(p)\land\\
\W(m,k,\bar x)\land \W(m,p,\bar x)\land \W(m,k,\bar y)\land
\W(m,p,y_1,y_2,y_3,t)\then y_4>t.
\end{multline*}
\end{thm}

\begin{figure}
\begin{center}
\psfrag*{text1}[c][c]{\shortstack[c]{my spaceship\\ is 1km long}}
\psfrag*{text2}[c][c]{\shortstack[c]{\ \ it's only $\sqrt{1-v^2}$\\
km long}} \psfrag*{text6}[t][t]{now ($m$)} \psfrag*{text7}[t][t]{1
second later ($m$)} \psfrag*{text5}[c][c]{\shortstack{1 km = 1
light-second\\ (in this picture)}}
\psfrag*{text3}[b][b]{$\sqrt{1-v^2}$}
\psfrag*{text4}[b][b]{$-v+\sqrt{1-v^2}$} \psfrag*{0}[b][b]{$0$}
\psfrag*{1}[l][l]{$1$} \psfrag*{m}[l][l]{\textcolor{gr}{$m$}}
\psfrag*{k}[l][l]{\AMbl{$k$}} \psfrag*{-v}[b][b]{$-v$}
\includegraphics[keepaspectratio, width=\textwidth]{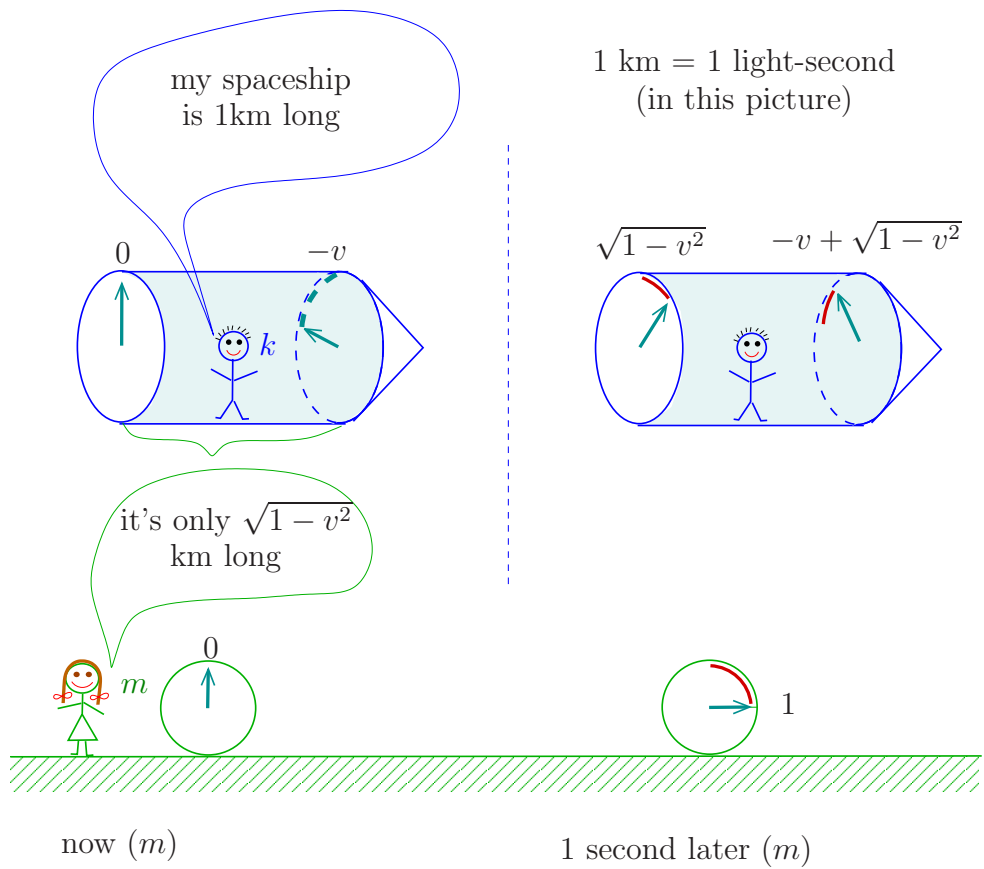}
\end{center}
\caption{\label{para1-fig}According to $m$, the length of the
spaceship is $\sqrt{1-v^2}$ km, it is $1$ km wide and tall, and the
clocks in the nose show $v$ less time than those in the rear, where
$v$ is the relative velocity of $m$ and $k$. According to $k$, the
length of the ship is $1$ km, it is 1 km wide and tall, and the
clocks in the nose and the ones in the rear all show the same time.}
\end{figure}

For {\bf proof} see, e.g., \cite[Thm.3.2.13]{Mphd}. What the average
layperson usually knows about the predictions of relativity is that
``moving clocks slow down, moving spaceships shrink, and moving
clocks get out of synchronism, i.e., the clock in the nose of a fast
moving spaceship is late (shows less time) when compared with the
clock in the rear." See Figure~\ref{para1-fig}. Let's call these
three predictions the ``paradigmatic effects" of SR. Now,
\ax{SpecRel} implies all the paradigmatic effects quantitatively,
too.\footnote{They follow from our next theorem. However, in our
works we usually prefer proving the paradigmatic effects one-by-one,
directly from the axioms of \ax{SpecRel} because this illuminates or
illustrates how we perform our conceptual analysis. These proofs can
be found in, e.g., \cite[sec.2.4]{logst}.} From this it follows that
the so-called worldview transformations are Poincar\'e-functions,
thus everything follows from our \ax{SpecRel} what follows from
``usual" special relativity theory.

Different observers may observe different spatial distance between
the same two events. This is so in Newtonian Kinematics (NK), too.
(For example, if I ate a sandwich and later drank a coffee on a train,
these two events were at the same place according to me, but
according to a coordinate system attached to Earth I ate the
sandwich at Budapest and drank the coffee at Vienna.) However, in NK
the time-difference between two events is the same for all
observers, it is ``absolute." According to the paradigmatic effects,
in RT even the time-difference between two events depends on the
state of motion of the observer! (The observer moving relative to
$m$ will observe less time passed between $e$ and $e'$ because his
clock ``slowed down.") In this respect, space and time in RT are
``more alike" than in NK. Our next theorem states that a certain
combination of spatial distance and time-difference is ``absolute"
in RT, too. %%
The proof of this theorem can be found in, e.g.,
\cite[p.650]{logst}. Let us define
\begin{equation*}
\mu(\bar x,\bar y) :=
(x_1-y_1)^2+(x_2-y_2)^2+(x_3-y_3)^2-(x_4-y_4)^2.
\end{equation*}
Thus $\mu(\bar x,\bar y)$ is squared spatial distance minus squared
time-difference between events $e$ and $e'$ if these events took
place at $\bar x$ and $\bar y$, respectively. This quantity is
called (squared) {\it relativistic distance} between events $e$ and
$e'$ (or squared Minkowski-distance between space-time locations
$\bar x$ and $\bar y$). According to the next theorem, relativistic
distance is ``absolute" in RT. In RT relativistic distance plays the
same role as absolute time in NK. Minkowski geometry is based on
relativistic distance (in place of Euclidean distance).

\begin{thm}
\begin{multline*}
\ax{SpecRel}\vdash
%\forall
%oo'x_1x'_1y_1y'_1z_1z'_1t_1t'_1x_2x'_2y_2y'_2z_2z'_2t_2t'_2\enskip
%\IOb(o)\land\IOb(o')\land \\
%\forall b\enskip (W(o,b,x_1,y_1,z_1,t_1)\leftrightarrow
%W(o',b,x'_1,y'_1,z'_1,t'_1))\land \\
%\forall b\enskip (W(o,b,x_2,y_2,z_2,t_2)\leftrightarrow
%W(o',b,x'_2,y'_2,z'_2,t'_2))\\ \then
%(x_1-x_2)^2+(y_1-y_2)^2+(z_1-z_2)^2-(t_1-t_2)^2\\ =
%(x'_1-x'_2)^2+(y'_1-y'_2)^2+(z'_1-z'_2)^2-(t'_1-t'_2)^2. $\\
%\smallskip
\forall oo'\bar x\bar x'\bar y\bar y'\enskip
\IOb(o)\land\IOb(o')\land \\
\forall b\enskip (\W(o,b,\bar x)\leftrightarrow \W(o',b,\bar
x'))\land
\forall b\enskip (\W(o,b,\bar y)\leftrightarrow \W(o',b,\bar y'))\\
%\then (x_1-y_1)^2+(x_2-y_2)^2+(x_3-y_3)^2-(x_4-y_4)^2\\ =
%(x'_1-y'_1)^2+(x'_2-y'_2)^2+(x'_3-y'_3)^2-(x'_4-y'_4)^2 .
\then \mu(\bar x,\bar y) = \mu(\bar x', \bar y').
\end{multline*}
\end{thm}

According to the theorem above, relativistic distance between events
is an ``absolute" property. Clearly, any property defined from it is
absolute, too. By the use of modern rigorous logic, it can be stated
and proved that the properties definable from relativistic distance
are the only absolute properties; and moreover all of \ax{SpecRel}
can be re-constructed from Minkowski geometry (i.e., from the
``pseudo-metric" $\mu$). How can one formulate such a statement
rigorously in formal logic?

Definability theory is one of the most beautiful parts of modern
logic, see \cite{chang-keisler}, \cite{Hodges}, \cite{Makkai},
\cite{Pambuccian}. It is about investigating connections between
theories formulated in FOL with completely different vocabularies
(such as, e.g., Finite Set Theory based on the $\epsilon$ relation
and Arithmetic based on $+$ and $\ast$). What would happen if we did
not consider, e.g., the quantities and operations on them as
primitives of the language? What happens if we are curious about
where these primitive notions come from, if we want to give them
``operational" meanings? What happens if we choose the so-called
causality relation as the only primitive symbol of our language (as,
e.g., in Robb \cite{Robb}, Mundy \cite{Mundy})? Can we compare then
these theories,  can we say that one of these is stronger or weaker
than the other, or that two such theories express the same amount of
``knowledge" about the world?

Definability theory is strongly related to relativity theory and to
positivists ideas. In fact, its existence was initiated by Hans
Reichenbach in 1924 \cite{Rei24}. Reichenbach in his works
emphasized the need of definability theory and made the first steps
in creating it. It was Alfred Tarski who later founded this branch
of mathematical logic. Since then it developed to a well-used and
powerful theory, in much extent due to the works of Michael Makkai.

Very briefly, the reason for the need of definability theory of
logic in relativity theory, as explained by Reichenbach, is as
follows. When one sets up a physical theory $Th$, one wants to use
only so-called observational\footnote{This observable/theoretical
hierarchy is not perfectly well defined and is known to be
problematic, but as Friedman \cite{friedman} puts it, it is still
better than nothing.} concepts, such as, e.g., ``meeting of two
particles." While investigating the (observationally based) theory
$Th$ (such as our $\ax{SpecRel}$), one defines new, so-called
``theoretical" concepts, such as, e.g., ``relativistic distance"
$\mu$. Some defined concepts then prove to be so useful that one
builds a new theory $Th'$ based on the most useful theoretical
concepts, and investigates this new theory $Th'$ in its own merits.
The new theory $Th'$ usually is simple, streamlined, elegant - built
in such a way that we satisfy our aesthetic desires. This is the
case with Minkowski geometry. The original theory $Th$ contains its
own interpretation, because we tried to use observational concepts.
The physical interpretation of the new streamlined theory $Th'$ is
its connection with $Th$. The strongest relationship between two
theories in different first-order languages is called definitional
equivalence. When two theories are definitionally equivalent, in the
rigorous sense of definability theory of FOL, the observational
oriented theory $Th$ can be recaptured completely from the
theoretical-oriented streamlined theory $Th'$ (and vice versa).

As examples, we can take our axiom system \ax{SpecRel} for $Th$, and we can take
Minkowski Geometry for $Th'$. Goldblatt \cite[App.A]{Goldblatt} gave
a complete FOL axiom system \ax{MG} for Minkowski geometry. His
theory is based on and is nicely analogous to Tarski's FOL
axiomatization for Euclidean geometry (see, e.g., \cite{TGgeom}).
The definitional equivalence of our present
``observationally-oriented" theory \ax{SpecRel} and the FOL theory
of Minkowski geometry \ax{MG} given
by Goldblatt is proved by Madar\'asz \cite[Chap.6.2]{Mphd}.%
\footnote{For proving this equivalence, one has to add extensionality
axioms for observers and light-signals to \ax{SpecRel}, and one has
to enrich \ax{MG} with a ``meter-rod."} It is noteworthy to mention
that in this application, relativity theory contributed to
definability theory once again: for the precise formulation of the
equivalence of the two theories we had to elaborate a methodology
for how to define new ``entities" (such as ``events") in addition to
the old methods which are about how to define new relations on
already existing entities (such as ``observer").
This definitional equivalence of the two theories can also be
expressed by saying that \ax{SpecRel} is complete with respect to
the Minkowskian model of SR generalized over ordered fields. Hence
everything which can be formulated in our language and true in these
Minkowskian models can be proved from our axiom system \ax{SpecRel}.
This is a kind of completeness theorem for the streamlined theory
\ax{SpecRel} with respect to Minkowskian Geometry as the intended
model for SR.

Less tight relationships than definitional equivalence between
theories are also very useful, these kinds of relationships are
called interpretability and duality connections. For an
illustration, let us turn to the question of where quantities and
coordinate systems come from. The axiom system \ax{AxSR} of James Ax
\cite{Ax} for SR is based on a first-order language that contains
only two unary relation symbols \mbox{\sf P,S} for ``particles" and
``signals" (corresponding to our ``bodies" and ``photons"), and two
binary relation symbols \mbox{\sf T,R} for ``transmitting a signal"
and
``receiving a signal.'' One can give an interpretation %\mbox{\sf int}
of our FOL theory \ax{SpecRel} in Ax's FOL theory \ax{AxSR} (see
\cite[proof-outline of Thm.2.1]{AMNsamples}). This amounts to
defining the primitive relations of the language of our \ax{SpecRel}
in terms of the primitives of Ax's \ax{AxSR}, and then proving from
\ax{AxSR} the translated axioms of \ax{SpecRel} as theorems. This is
an interpretation in the sense of definability theory. Now, this
interpretation also can be thought of as giving a kind of
operational ``definition" for how to set up ``operationally" the
coordinate systems appearing in \ax{SpecRel} as primitives. The
question of how to give algorithms for setting up coordinate systems
in this context is treated in more detail and depth in Szab\'o
\cite{Szabo}.

Theories form a rich structure when we investigate their
interconnections.  G\"odel's incompleteness theorem pointed already
in the direction of investigating hierarchies of theories rather
than single theories. (There is no ``strongest" theory for the
interesting subjects, there are only stronger and stronger
theories.) Answering ``why-questions," ``reverse mathematics,"
modularizing our knowledge all point to the study of weaker and
weaker theories, and also to studying the interpretations between
theories (see \cite{myphd}). Algebraic logic, developed by Tarski
and his followers, is a branch of definability theory which
establishes a duality between hierarchies of theories and between
classes of algebras (cf., e.g., \cite[Chap.4.3]{HMTII},
\cite{HMTAN}, \cite{Makkai}). In modern approaches to logic,
theories are considered as dynamic objects as opposed to the more
traditional ``eternally frozen" idea of theories. For approaches to
the dynamic trend in mathematical logic cf., e.g., van Benthem
\cite{vB96}, Gabbay \cite{Gabbay}, and \cite{MNS}. This new
``plurality of theories" or ``hierarchy of small theories (as
opposed to a single monolithic one)" approach can help realizing the
central or essential VC-aims without the old stumble blocks of the
original VC attempts. This is a wisdom gained from FOM, see
\cite{FOM}.

In the rest of the paper we briefly indicate how to arrive to a
transparent FOL axiomatization of general relativity from our
\ax{SpecRel}. By this we realize Einstein's original program
formally and literally.

\section{first step toward gr: extending the theory to accelerated observers}

\comment{There are quite a few axiomatizations of Special Relativity
(SR) in the literature, see, e.g., \cite{Guts}, \cite{Mundy},
\cite{Robb}, \cite{Schutz73}, \cite{Schutz}, \cite{Suppes}. Two of
them are formulated within FOL, namely \cite{Ax} and
\cite{Goldblatt}, \cite{Pambuccian}... However, there are only a few
axiomatizations of General Relativity (GR), see, e.g.,
\cite{Mould}..., and a very few of them (besides ours) are
formulated within FOL, see, e.g., \cite{Basri}, \cite{Benda}.
Moreover, as far as we know, no paper in the literature deals with
relations of the axioms of the two theories. Here we are going to
fill this gap by deriving the axioms of GR of \cite{logst},
\cite{myphd} from the axioms of SR of \cite{AMNsamples} in two
simple steps.  Both of these steps are natural and well motivated by
some central ideas of Einstein's. In the first step we are going to
extend our axiomatization of SR to accelerated observers. This step
provides us a theory which nicely fills the gap between SR and GR.
In the second step we are going to arrive to a first-order theory of
GR by eliminating the distinction between inertial and non-inertial
observers in the level of axioms.}%endcomment

As a first step toward General Relativity theory (GR), we are going
to extend our \ax{SpecRel} theory with accelerated observers.  By
accelerated observer we mean any not necessarily inertial observer.
Let us first note that none of the axioms of \ax{SpecRel} speaks
about noninertial observers.

Since in the language we have already introduced the concept of
arbitrary observer, the only thing we have to do is to assume some
axioms about them. Our key axiom to assume about arbitrary observers
is the following:
\begin{description}
\item[\ax{AxCmv}:] At each moment of its life, any observer
  coordinatizes (``sees") the nearby world for a short while in the same way as
  some inertial observer does.
\end{description}
For precise formulation of this axiom in the spirit we formulated
the axioms of \ax{SpecRel} see \cite{logst}, \cite{twp},
\cite{myphd}. Let \ax{AccRel}$^-$ be the axiom system consisting of
\ax{AxCmv} and all the axioms of \ax{SpecRel}.

Let us see how strong our theory \ax{AccRel}$^-$ is. To test its
strength we are going to investigate whether the Twin Paradox (TwP)
and the gravitational time dilation (``gravity causing slow time")
are provable from it.

According to TwP, if a twin makes a journey into space
(accelerates), he will return to find that he has aged less than his
twin brother who stayed at home (did not accelerate). However much
surprising TwP is, it is not a contradiction. It is only a fact
showing that the concept of time is not as simple as it seems at
first.

A more optimistic consequence of Twp is the following. Suppose you
would like to visit a distant galaxy 200 light years away. You are
told it is impossible because even light travels there for 200
years. But you do not despair, you accelerate your spaceship nearly
to the speed of light. Then you travel there in 1 year subjective
time. When you arrive back, you aged only 2 years. So you are happy,
but of course you cannot tell the story to your brother, who stayed
on Earth. Alas you can tell it to your grand-\ldots-grand-children
only.

In the FOL language introduced in this paper we can formulate TwP,
see \cite{twp}, \cite{myphd}. Let us denote the formulated version
of TwP as \ax{TwP}.

\ax{AccRel}$^-$ is not yet strong enough to imply \ax{TwP}. One
would think that this is so because we did not state enough
properties of the real numbers for speaking about curved lines.
However, even assuming \ax{Th(\R)}, i.e., all the FOL formulas valid
in the real numbers, together with \ax{AccRel}$^-$ is not sufficient
to prove
\ax{TwP}, %i.e., the following can be proved,
see \cite{twp}, \cite{myphd}:
\begin{thm}
$\ax{AccRel}^-\cup\ax{Th(\R)}\nvdash \ax{TwP}.$
\end{thm}

We note that the above theorem is a theorem stating that one cannot
prove \ax{TwP} from $\ax{AccRel}^-\cup\ax{Th(\R)}$, it is not only
the case that we are not ``clever enough" to find a proof but there
is none. Its proof goes via using the completeness theorem of FOL,
namely we find a model in which all the formulas in
$\ax{AccRel}^-\cup\ax{Th(\R)}$ are true, but in which \ax{TwP} is
not true.

This theorem states that even assuming every first-order formula
which is true in $\R$ is not enough for our purposes. At first sight
this result suggests that our programme of FOL axiomatization of GR
breaks down at the level of TwP. It would be depressing if we were
not able to keep our axiomatization within FOL, because there are
weighty methodological reasons for staying within it, see, e.g.,
\cite[Appendix]{pezsgo}, \cite[sec.11]{myphd}.  However, we are
saved: in our language there is a FOL axiom scheme (nice set of
axioms) called \ax{IND} which is sufficient for our purposes. Axiom
scheme \ax{IND} expresses that every nonempty and bounded subset of
the quantities which is parametrically definable in our language has
a least upper bound (i.e., supremum). \ax{IND} is a first-order
logic approximation of the second-order logic continuity axiom of
the real numbers, and it belongs to the methodology developed in FOM
and in reverse mathematics that \ax{AxField} strengthened with
\ax{IND} are strong enough for a FOL treatment of areas involving
the real numbers.

Together with this scheme \ax{AccRel}$^-$ implies \ax{TwP}, i.e.,
the following theorem can be proved, see \cite{twp}, \cite{myphd}:
\begin{thm}$\ax{AccRel}^-\cup\ax{IND}\vdash \ax{TwP}.$
\end{thm}

How can a FOL axiom scheme be stronger than all the FOL formulas
valid in $\R$?  The answer is that \ax{IND} is formulated in a
richer language than that of the reals, hence it can state more than
the whole FOL theory of $\R$. If we assume \ax{IND} only for
formulas in the language of ordered fields, we get an axiom schema
equivalent to \ax{Th(\R)}, see \cite{myphd}.\footnote{Actually, the
restriction of \ax{IND} to fields $\langle Q,+,\ast,<\rangle$
coincides with Tarski's FOL version of Hilbert's continuity axiom
for geometry, cf.\ \cite[p.71, axiom B5]{Goldblatt}.} Let us now
introduce our axiom system for accelerated observers as:
\begin{equation*}
\boxed{\ax{AccRel}=\ax{SpecRel}\cup\{\ax{AxCmv}\}\cup\ax{IND}}.
\end{equation*}

Let us continue with the gravitational time dilation. By Einstein's
equivalence principle, we can also formulate the statement ``gravity
causes slow time'' (usually called ``gravitational time dilation"
GTD) in our language. Moreover, the formulated version of this
statement is provable from the theory $\ax{AccRel}$, see
\cite{MNSz}, \cite{myphd}. The \ax{AccRel} formulation of GTD
basically says that in any accelerated spaceship the clocks in the
rear run slower than those in the nose. (The effect is increasing
with increasing acceleration. Moreover, it approaches infinity as
the acceleration does.) So we are able to derive nontrivial
predictions about gravity before we have introduced any axiom system
of GR.

The theory \ax{AccRel} is halfway between SR and GR. Einstein used a
non-formalized version of \ax{AccRel} as a heuristic in introducing
GR, e.g., when he made predictions about the influence of
gravitation on the propagation of light \cite{E11}, \cite[\S\S
18-22]{E21}.

\section{second step: ``emancipating'' noninertial observers}

We are going to modify the axioms of \ax{SpecRel} and \ax{AxCmv} one
by one and get an axiomatic theory of general relativity. The
modification consists of ``eliminating the privileged class of
inertial reference frames,'' which was a central idea of Einstein's,
see \cite[\S\S 18-22]{E21}, \cite{friedman}. We replace each axiom
of \ax{SpecRel} by a new one which does not speak about inertiality
but otherwise the content of which tries to approximate that of the
old one. All the new axioms will be motivated by our theory
\ax{AccRel}. Roughly, each axiom of \ax{AccRel} will be replaced by
a ``generalized" version which does not mention inertiality and
which is still in the spirit of \ax{AccRel}.

The generalized version of \ax{AxSelf} is the following:
\begin{description}
\item[\ax{AxSelf^-}:] An observer coordinatizes itself on a subset of the time axis:
\begin{equation*}
\forall oxzyt\enskip \W(o,o,x,y,z,t)\then x=y=z=0.
\end{equation*}
\end{description}

The modified version \ax{AxEv}$^-$ of \ax{AxEv} contains the
following two statements: (1) any observer coordinatizes the events
in which it was observed by some other observer, and (2) if observer
$o$ coordinatizes an event which is coordinatized by observer $o'$,
then $o$ also coordinatizes the events which are near this event
according to $o'$. This can be summarized as follows:
\begin{description}
\item[\ax{AxEv^-}:] Any observer coordinatizes the events in which it was
 observed; and the domains of worldview transformations are open.
\end{description}
The modified versions of \ax{AxPh} and \ax{AxSymd} are achieved by
localizing and generalizing them, i.e., we get the modified versions
by restating these axioms only in infinitesimally small
neighborhoods, but for every observer. The idea that ``GR is locally
SR'' also goes back to Einstein. Our symmetry axiom \ax{AxSymd} has
many equivalent versions, see \cite[sec.s 2.8, 3.9, 4.2]{pezsgo}. We
can localize any of these versions and use it in a FOL axiom system
for GR. For aesthetic reasons here we localize \ax{AxSymt}, the
version stating that inertial observers see each others' clocks
behave the same way. So \ax{AxPh^-} and \ax{AxSymt^-} are the
formalized versions of the following statements:
\begin{description}
\item[\ax{AxPh^-}:] The instantaneous velocity of photons is $1$ in
 the moment when they ``meet" the observer who coordinatizes
 them, and any observer can send out photons in any direction
 with this instantaneous velocity.
\end{description}

\begin{description}
  \item[\ax{AxSymt^-}:] Meeting observers see each other's clocks
  behaving the same way, %slowing down with the same ratio,
  at the event of meeting.
\end{description}
For formulation of these axioms and the corresponding concepts
in our first-order language, see \cite{myphd}.

Now all the four axioms of \ax{SpecRel} are modified according to
the above requirements.  Strictly following these guidelines,
\ax{AxCmv^-} would state that the worldview transformation between
observers are differentiable in their meeting-point. To avoid
baroque, we state simply differentiability of the worldview
transformations. A natural generalization is $n$-times
differentiability (which is natural to consider in view of our
wanting to speak about location, speed and acceleration). Each axiom
of this series of potential axioms can be formulated in the language
above by the techniques used in \cite{logst}, \cite{twp},
\cite{myphd}.
\begin{description}
\item[\ax{AxDiff_n}:] The worldview transformations are $n$-times
 differentiable functions.
\end{description}
Let us introduce the following simple axiom systems for general
relativity:
\begin{equation*}\index{\ax{GenRel_n}}
\boxed{\ax{GenRel_n}\leteq \Setopen \ax{AxField}, \ax{AxSelf^-},
\ax{AxPh^-}, \ax{AxEv^-},\ax{AxSymt^-},\ax{AxDiff_n}
\Setclose\cup\ax{IND}}.%\cup\ax{IND}.
\end{equation*}
%Let us note that any model of \ax{GenRel_m} is a model of
%\ax{GenRel_n} if $m\ge n$.
The following theorem illustrates that our axiom system
\ax{GenRel_n} captures the $n$-times differentiable standard models
of usual GR well.

Lorentzian manifolds are the intended models of GR, much the same
way as Minkowski geometry was the intended model of SR. Roughly, a
Lorentzian manifold is a geometry which at every of its points
locally looks like the Minkowski geometry, cf., e.g.,
\cite[p.23]{Wald}.

\begin{thm}\label{thm-gr}
\ax{GenRel_n} is complete with respect to the $n$-times
differentiable Lorentzian manifolds over real-closed fields.
\end{thm}

\comment{Let us also introduce an infinite version of \ax{GenRel}
which contains all the axioms of the series above:
\begin{equation*}\index{\ax{GenRel_\infty}}
%\boxed{
\ax{GenRel_\infty}\leteq \bigcup_{n\ge 1}\ax{GenRel_n}%}
\end{equation*}
By the following theorem \ax{GenRel_\infty} captures the smooth models of GR well.
\begin{thm}\label{thm-grinf}
\ax{GenRel_\infty} is complete with respect to smooth Lorentzian
manifolds over ordered fields.
\end{thm}}

There are many interesting GR space-times, black holes, worm-holes,
time-warps, etc.  The physical relevance of these so called exotic
space-times increases with time. For instance, there is a rapidly
growing number of experimental evidence for huge slowly rotating
black holes, which are the
simplest examples of time-warps. By Theorem \ref{thm-gr} %and
%\ref{thm-grinf}
even the most exotic model of GR is also a model of
our \ax{GenRel} theory. Hence, within \ax{GenRel} we can investigate
the properties of these exotic models.

To ensure that we can do indeed physics in the framework of
\ax{GenRel}$_n$ ($n\ge 3$) we defined in \cite{logst}, \cite{myphd}
the notion of time-like geodesics in terms of \ax{GenRel}. These
serve as world-lines of inertial bodies. So, though we abandoned
inertial observers as primitives, inertial motion becomes
accessible/definable as a derived notion (in terms of the primitives
of \ax{GenRel}). Space-time curvature is defined from geodesics the
usual way. So, in particular, the outcomes of experiments involving
inertial motion can be predicted (e.g., computing the trajectories
of bullets or photon geodesics) on the basics of the new,
streamlined theory \ax{GenRel} in a purely logical way.

\section{concluding remarks} As it was the case with \ax{SpecRel},
cf., \cite{AAMN}-\cite{dyn-studia}, having obtained the streamlined
axiomatization \ax{GenRel} and its completeness for ``usual" GR is
only a first step towards a logic based conceptual analysis of GR,
its predictions, alternatives or variants, answering the
why-questions in a spirit which is a natural continuation of the VC
programme.

\section*{Acknowledgements} We thank the organizers and the
participants of the workshop ``Vienna Circle and Hungary" held in
Vienna in March 2008, for their hospitality, questions, and remarks.
We also thank Rainer Tiemeyer for inspiring letters about
formalizing scientific theories. Research supported by the National
Foundation for Scientific Research grant No T73601 as well as by a
Bolyai Grant for Judit X.\ Madar\'asz.
% \eject

%\theendnotes

%\section{trash}
Alfr\'ed R\'enyi Institute of Mathematics\\
of the Hungarian Academy of Sciences\\
Budapest P.O. Box 127, H-1364 Hungary.\\
andreka@renyi.hu, madarasz@renyi.hu, nemeti@renyi.hu,
nemetip@gmail.com, turms@renyi.hu


\begin{thebibliography}{99}
\footnotesize

\bibitem{Ax} Ax, J., The elementary foundations of spacetime,
  {\it Foundations of Physics} 8, 7-8 (1978), pp.507-546.

\bibitem{AAMN} Andai, A., Andr{\'e}ka, H., Madar{\'a}sz, J. X.\ and
N{\'e}meti, I., Visualizing some ideas  about G\"odel-type rotating
universes, in: M. Scherfner and M. Plaue, eds., {\it G\"odel-type
spacetimes: history and new developments}, Springer, to appear.
arXiv:0811.2910v1.

\bibitem{AMNsamples} Andr{\'e}ka, H., Madar{\'a}sz, J. X.\ and
  N{\'e}meti, I., Logical axiomatizations of space-time. Samples
  from the literature, in: A. Pr\'ekopa and E. Moln\'ar, eds., {\it Non-Euclidean
    geometries}, volume 581 of {\it Mathematics and Its
    Applications}. Springer, 2006, pp.155-185. %?OK??

\bibitem{logst} Andr{\'e}ka, H., Madar{\'a}sz, J. X.\ and
  N{\'e}meti, I., Logic of space-time and relativity
  theory, in: M. Aiello, I. Pratt-Hartmann and J. van Benthem,
  eds., {\it Handbook of Spatial Logics}. Springer, 2007,
  pp.607-711.

\bibitem{pezsgo} Andr{\'e}ka, H., Madar{\'a}sz, J. X.\ and
  N{\'e}meti, I., {\it On the logical structure of relativity
    theories}.  E-book, Alfr{\'e}d R{\'e}nyi Institute of
  Mathematics, Budapest, 2002.  With contributions from Andai, A.,
  S{\'a}gi, G., Sain, I.\ and T{\H o}ke, Cs.
  http://www.math-inst.hu/pub/algebraic-logic/Contents.html. 1312
  pp.

\bibitem{dyn-studia} Andr{\'e}ka, H., Madar{\'a}sz, J. X.,
  N{\'e}meti, I.\ and Sz{\'e}kely, G., Axiomatizing
  relativistic dynamics without conservation postulates, {\it
    Studia Logica} 89, 2 (2008), pp.163-186.

%    \bibitem{ANS01} Andreka, H., N\'emeti, I., and Sain, I.,
%    Algebraic Logic. In: ...

\bibitem{ANW} Andr{\'e}ka, H., N{\'e}meti, I.\ and W{\"u}thrich, C.,
A twist in the geometry of rotating black holes:
  seeking the cause of acausality,
 {\it General Relativity and Gravitation} 40, 9 (2008), pp.1809-1823.

%\bibitem{Basri} Saul A. Basri, {\it A Deductive Theory of Space and
%  Time}. Amsterdam: North-Holland 1966.

%\bibitem{Benda} Thomas Benda, ``A Formal Construction of the Spacetime
%  Manifold'', in: {\it Journal of Philosophical Logic} 37, 5, 2008,
%  pp.441-478.

  \bibitem{uni1} Carnap, R.,  Frank, P., Hahn, H., Neurath, O., Joergensen, J.\ and Morris, C. (eds).,
  {\it Unified Science}.
  These works are translated in {\it Unified Science: The Vienna Circle Monograph Series
  Originally Edited by Otto Neurath}, Kluwer, 1987.

\bibitem{chang-keisler} Chang, C. C.\ and H. Keisler, J., {\it
  Model theory}. North-Holland 1973, 1977, 1990.

\bibitem{E11} Einstein, A., \"Uber den Einfluss der Schwercraft auf
  die Ausbreitung des Lichtes, {\it Annalen der Physik} 35, 10 (1911), pp.898-908.

  \bibitem{E21} Einstein, A., \"Uber die Spezielle und die Allgemeine
  Relativit{\"a}tstheorie, {\it Verlag von F. Vieweg and Son, Braunschweig}, 1921.

%\bibitem{Etesi-Nemeti} G{\'a}bor Etesi and Istv{\'a}n N{\'e}meti,
%  ``Non-Turing computations via Malament-Hogarth space-times'', in: {\it
%  International Journal of Theoretical Physics} 41, 2, 2002,
%  pp.341-370.

 \bibitem{FOM} Friedman, H., On foundational thinking 1, {\it Posting
in FOM (Foundations of Mathematics)}, Archives www.cs.nyu.edu,
January 20 2004.

\bibitem{friedman} Friedman, M., {\it Foundations of Space-Time
  Theories. Relativistic Physics and Philosophy of Science},
  Princeton University Press 1983.

  \bibitem{friedman2} Friedman, M., {\it Reconsidering logical
  positivism},  Cambridge University Press 1999.

  \bibitem{Gabbay} Gabbay, D. M., {\it Labelled Deductive Systems},
  Clarendon Press, 1996.

\bibitem{Goldblatt} Goldblatt, R., {\it Orthogonality and spacetime
  geometry}, Springer-Verlag, 1987.

  \bibitem{godel} G\"odel, K., An Example of a New Type of
  Cosmological Solutions of Einstein's Field Equations of
  Gravitation, {\it Reviews of Modern Physics} 21, 3 (1949),
  pp.447-450.

%\bibitem{Guts} Alexander K.\ Guts, ``The axiomatic theory of relativity'', in:
%{\it Russian Mathematical Surveys} 37, 2, 1982, pp.41-89.

%\bibitem{Gyenis-Roberts} Bal{\'a}zs Gyenis and Bryan W. Roberts ``Supertasks:
%  G{\"o}del strikes back'', Manuscript, University of Pittsburgh, 2007.\?

\bibitem{HMTII} Henkin, L., Monk, J. D., and Tarski, A., {\it Cylindric
Algebras Part II}, North-Holland, 1985.

\bibitem{HMTAN} Henkin, L., Monk, J. D., Tarski, A., Andr\'eka, H.,
and N\'emeti, I., {\it Cylindric Set Algebras}, Lecture Notes in
Mathematics Vol 883, Springer Verlag, 1981.

\bibitem{Hirsch} Hirsch, R., Einstein: Logic, Philosophy,
Politics. Abstract for talk in ``Logic in Hungary, 2005", Budapest,
Hungary August 5-10, 2005.
http://atlas-conferences.com/cgi-bin/abstract/caqb-48

\bibitem{Hodges} Hodges, W., {\it Model Theory}, Cambridge University Press, 1993.

%  \bibitem{Ho05} Horv\'ath, R., An Alexandrov-Zeeman tpe theorem and
%  relativity theory. Paper for Scientific Student Contest, E
%  "otv\"os Lor\'and University, Budapest, 2005.

\bibitem{Mphd} Madar{\'a}sz, J. X., {\em Logic and Relativity (in
  the light of definability theory)}.  PhD thesis, E{\"o}tv{\"o}s
  Lor{\'a}nd Univ., Budapest, 2002.
  http://www.math-inst.hu/pub/algebraic-logic/diszi0226.pdf.gz

\bibitem{twp} Madar{\'a}sz, J. X., N{\'e}meti, I.\ and
  Sz{\'e}kely, G., Twin paradox and the logical foundation of
  relativity theory, {\it Foundations of Physics} 36, 5 (2006),
  pp.681-714.

\bibitem{MNSz} Madar{\'a}sz, J. X., N{\'e}meti, I.\ and Sz{\'e}kely, G.,
First-order logic foundation of relativity theories, in:
  D. M. Gabbay, S. Goncharov and M. Zakharyaschev, eds., {\it
    Mathematical problems from applied logic II}. Springer,
  2007, pp.217-252.

%\bibitem{MNT} Judit X.~Madar{\'a}sz, Istv\'an N{\'e}meti, and Csaba
%  T{\H o}ke, ``On generalizing the logic-approach to space-time
 % towards general relativity: first steps'', in: Vincent F.\ Hendricks,
%  Fabian Neuhaus, Stig Andur Pedersen, Uwe Scheffler and Heinrich Wansing
%  (Eds.), {\it First-Order Logic Revisited}.  Berlin: Logos Verlag
%  2004, pp.225-268.

  \bibitem{Makkai} Makkai, M., {\it Duality and definability in first
  order logic}, Memoirs of the American Mathematical Society No 503,
  Providence, Rhode Island, 1993.

%\bibitem{Mould} Richard A. Mould, ``An Axiomatization of Generla
%  Relativity'', in: {\it Proceedings of the American Philosophical
%    Society} 103, 2, 1959, pp.485-529.

\bibitem{MNS} Mikul\'as, Sz., N\'emeti, I. and Sain, I., Decidable
Logics of the Dynamic Trend, and Relativized Relation Algebras. In:
L. Csirmaz, D. M. Gabbay and M. de Rijke, eds., {\it Logic
Colloquium'92}, Studies in Logic, Language and Computation, CSLI
Publications, 1995. pp.165-175.

\bibitem{Mundy} Mundy, B., Optical axiomatization of {M}inkowski
  space-time geometry, {\it Philosophy of Science} 53, 1 (1986),
  pp.1-30.

\bibitem{Nemeti-dgy} N{\'e}meti, I. and D{\'a}vid, Gy.,
  ``Relativistic computers and the Turing barrier'', in: {\it Applied
  Mathematics and Computation} 178, 1 (2006), pp.118-142.

   \bibitem{uni2} Neurath, O., Carnap, R.\ and Morris, C. (eds),
  {\it Foundations of the Unity of Science}. The university of Chicago Press,
  Series: Foundations of the Unity of  Science: Toward an International
  Encyclopedia of Unified Science. Vol.s\ 1 and 2.\ 1971.

\bibitem{Pambuccian} Pambuccian, V., Alexandrov-{Z}eeman type
  theorems expressed in terms of definability, {\it Aequationes
    Mathematicae} 74, 3 (2007), pp.249-261.

\bibitem{Rei24} Reichenbach, M., {\it Axiomatization of the theory of relativity},
Univeversity of California Press, Berkeley, 1969.
Translated by M. Reichenbach. Original German edition published in 1924.

\bibitem{Robb} Robb, A. A., {\it A Theory of Time and
  Space}. Cambridge University Press, 1914.

%\bibitem{Schutz73} John W.\ Schutz, {\it Foundations of special
%  relativity: kinematic axioms for {M}inkowski space-time} Berlin: Springer-Verlag 1973.

%\bibitem{Schutz} John W.\ Schutz, ``An axiomatic system for
%  {M}inkowski space-time'', in: {\it Journal of Mathematical Physics}
%  22, 2, 1981, pp.293-302.

%\bibitem{Suppes}  Patrick Suppes, ``Axioms for relativistic kinematics
%with or without parity'', in: Leon Henkin, Patrick Suppes and Alfred Tarski (Eds.),
%{\it The axiomatic method. With special reference to geometry and physics}.
%Berkeley: Nort-Holland 1959, pp.291-307.

\bibitem{Szabo} Szab\'o, L. E., Empirical foundation of space and time.
In: M. Su\'arez, M. Dorato and M. R\'edei, eds., {\it  EPSA;
Epistemology and Methodology: Launch of the European Philosophy of
Science Association Proceedings of the First Conference of the
European Philosophy of Science Association, Madrid, Spain, November
14-17, 2007}, Springer, 2009.

\bibitem{Szekely} Sz\'ekely, G.,  On why-questions in physics,
in this volume.

\bibitem{myphd} Sz{\'e}kely, G., {\it First-Order Logic
  Investigation of Relativity Theory with an Emphasis on Accelerated
  Observers}, PhD thesis, E{\"o}tv{\"o}s Lor{\'a}nd Univ., Budapest,
  2009. %?OK??

  \bibitem{Ta59} Tarski, A., What is elementary geometry?, In:
  L. Henkin, P. Suppes and A. Tarski, eds.,  {\it The
  axiomatic method with special reference to geometry and physics},
  North-Holland, 1959, pp.16-29.

  \bibitem{TGgeom} Tarski, A., and Givant, S. R., Tarski's system of
  geometry, {\it Bulletin of Symbolic Logic}, 5, 2 (1999), pp.175-214.

  \bibitem{vB96} van Benthem, J. F. A. K., {\it Exploring logical
  dynamics}, Sudies in Logic, Language and Information, CSLI
  Publications, 1996.

  \bibitem{Wald} Wald, R. M., {\it General Relativity}, The
  University of Chicago Press, 1984.

\end{thebibliography}
\end{document}